\documentstyle[twoside,12pt]{article} 
\def\N{I\!\!N}

\def\cR{{\cal R}} 
  
\newtheorem{theorem}{Theorem} 
\newtheorem{lemma}{Lemma} 
\newtheorem{corollary}{Corollary} 
 
\newtheorem{definition}{Definition} 
 
\newtheorem{remark}{Remark} 
\newtheorem{question}{Question}

\newcommand{\mod}{{\rm \; mod\ }} 
 
\date{} 

\title{Reidemeister number of any automorphism  of a Gromov hyperbolic group is infinite}
\author{Alexander Fel'shtyn\thanks{Part of this work was conducted during  author  stay in Universite Paul Sabatier, Toulouse  and Ohio State University, Columbus}}
\begin{document}   
\bibliography{ref} 
\bibliographystyle{plain}
\maketitle

 \bigskip
\begin{abstract}  
We show that the number of twisted conjugacy classes is infinite for any automorphism 
of non-elementary, Gromov hyperbolic group \cite{g}. An analog of Selberg theory 
for twisted conjugacy classes is proposed. 

\end{abstract}

\setcounter{section}{-1}
 
\section{Introduction}

Let $G$ be a finitely generated group and $\phi: G\rightarrow G$ an endomorphism.
Two elements $\alpha,\alpha^\prime\in G$ are said to be
$\phi-conjugate$ iff there exists $\gamma \in G$ with
$$
\alpha^\prime=\gamma  \alpha   \phi(\gamma)^{-1}.
$$
We shall write $\{x\}_\phi$ for the $\phi$-conjugacy class
 of the element $x\in G$.
The number of $\phi$-conjugacy classes is called the $Reidemeister$
$number$ of an  endomorphism $\phi$, denoted by $R(\phi)$. If $\phi$ is the identity map then the $\phi$-conjugacy classes are the usual conjugacy classes in the group $G$. 

We note that $R(\phi)$ is infinite if group $G$ is free Abelian and
the action of $\phi$ on $G$ has $1$ as eigenvalue  \cite{f1} .
 
In \cite{fh} we have conjectured that the Reidemeister number is   infinite
as long as the endomorphism $\phi$ is injective and  the group $G$ has exponential growth.

In this paper  we prove our conjecture for  any automorphism of any    non-elementary
(i.e. not virtually cyclic), Gromov hyperbolic group.
We also prove some generalisations of this result. 

Main result of this paper has  topological counterpart.
 Let $X$ to be a connected, compact
polyhedron and $f:X\rightarrow X$ to be a continuous map.
Let $p:\tilde{X}\rightarrow X$ be the universal cover of $X$
and $\tilde{f}:\tilde{X}\rightarrow \tilde{X}$ a lifting
of $f$, i.e. $p\circ\tilde{f}=f\circ p$.
Two liftings $\tilde{f}$ and $\tilde{f}^\prime$ are called
{\sl conjugate} if there is a element$\gamma$ in the deck transformation group $\Gamma\cong\pi_1(X)$
such that $\tilde{f}^\prime = \gamma\circ\tilde{f}\circ\gamma^{-1}$.
The subset $p(Fix(\tilde{f}))\subset Fix(f)$ is called
{\sl the fixed point class of $f$ determined by the lifting class $[\tilde{f}]$}.
Two fixed points $x_0$ and $x_1$ of $f$ belong to the same fixed point class iff
 there is a path $c$ from $x_0$ to $x_1$ such that $c \cong f\circ c $ (homotopy relative endpoints). This fact can be considered as an equivalent definition of a non-empty fixed point class.
 Every map $f$  has only finitely many non-empty fixed point classes, each a compact
 subset of $X$. 

A fixed point class is called {\sl essential} if its index is nonzero.
The number of lifting classes of $f$ (and hence the number
of fixed point classes, empty or not) is called the {\sl Reidemeister number} of $f$,
denoted $R(f)$.
This is a positive integer or infinity.
The number of essential fixed point classes is called the {\sl Nielsen number}
of $f$, denoted by $N(f)$.The Nielsen number is always finite.

 Theorem 3 implies that the topological Reidemeister number $R(f)$ is infinite for  a homeomorphism $f$
of a compact polyhedron $X$ with a non-elementary, Gromov hyperbolic fundamental group $\pi_1(X,x_0)$.

In the case of pseudo-Anosov homeomorphisms of surfaces we can develop
an analog of the Selberg theory and we obtain an asymptotic  expansion for 
the number of twisted conjugacy classes or  for the number of 
Nielsen fixed point classes whose norm is at most $x$.

The author would like to thank M. Bridson, B. Bowditch, D. Burghelea, M. Davis,  R. Grigorchyk, M. Gromov, G. Levitt, M. Lustig, Ch. Epstein, R. Hill, L. Potyagailo,  for stimulating discussions and comments.
The author would like to thank Laboratoire de Matematiques Emille Picard, Toulouse and the Departments of Mathematics, Ohio State University for their kind hospitality and support.

\section{Twisted conjugacy classes and Reidemeister number of group endomorphism}
   
\begin{lemma}
If $G$ is a group and $\phi$
 is an endomorphism of $G$ then
 an element $x\in G$ is always
 $\phi$-conjugate to its image $\phi(x)$.
\end{lemma}
{\sc  Proof.}
If $\gamma=x^{-1}$,  then one has immediately
 $\gamma x = \phi(x) \phi(\gamma)$.
The existence of a $\gamma$ satisfying this equation implies
 that $x$ and $\phi(x)$ are $\phi$-conjugate.

The mapping torus $M(\phi)$ of the group endomorphism $\phi: G\rightarrow G$
is obtained from group $G$ by adding a new generator $z$  and
adding the relations  $z^{-1}gz=\phi(g)$ for all $g\in G$. 
This means that $M(\phi)$ is a semi-direct product of $G$ with $Z$.

\begin{lemma}
Two elements $x,y$ of $G$ are $\phi$-conjugate if and only if  $xz$ and $yz$ are conjugate 
in the usual sense in $M(\phi)$. Therefore $R(\phi)$ is the number of usual conjugacy 
classes in the coset $G\cdot z$ of $G$ in $M(\phi)$. 
 \end{lemma}
{\sc  Proof.}
If $x$ and $y$ are $\phi$-conjugate,  then there is a $\gamma \in G$ such 
that $\gamma x=y\phi(\gamma)$. This implies $\gamma x=yz\gamma z^{-1}$ 
and therefore $\gamma(xz)=(yz)\gamma$.   So $xz$ and $yz$ are conjugate in the usual sense in $M(\phi)$. 
Conversely suppose $xz$ and $yz$ are conjugate in $M(\phi)$. 
Then there is a $\gamma z^n \in M(\phi)$ with $\gamma z^n xz=yz\gamma z^n$. 
From the relation $zxz^{-1}=\phi(x)$,  we 
obtain $\gamma \phi^n(x)z^{n+1}=y\phi(\gamma) z^{n+1} $ 
 and therefore $\gamma \phi^n(x)=y\phi(\gamma)$. 
This shows that $\phi^n(x)$ and $y$ are $\phi$-conjugate. 
However,  by lemma 1 ,   $x$ and $\phi^n(x)$ are $\phi$-conjugate, so $x$ and $y$ 
must be $\phi$-conjugate.

\begin{lemma}(T.Delzant)
Let $J$ be a non-elementary, Gromov hyperbolic group. Let $K$ be a normal
subgroup with abelian quotient. Then every coset $C$ of $J  \mod  K$ contains 
infinitely many conjugacy classes.
\end{lemma}
{\sc Proof.} (See \cite{LL}).
Fix $u$ in the coset $C$ under consideration. Suppose for a moment that
we can find $c,d \in K$, generating a free group of rank 2, such that
$uc^{\infty}\not=c^{-\infty}$ and $ud^{\infty}\not=d^{-\infty}$
(recall that we denote $g^{-\infty}=\lim_{n\rightarrow +\infty}g^{-n}$
for $g$ of infinite order). Consider $x_k=c^kuc^k$ and $y_k=d^kud^k$.
For $k$ large, the above inequalities imply that these two elements have infinite order, and do not generate a virtually cyclic group because ${x_k}^{+\infty}$ and ${x_k}^{-\infty}$(respectively ${y_k}^{+\infty}$ and ${y_k}^{-\infty}$ ) is close  to $c^{+\infty}$ and $c^{-\infty}$ (
respectively $d^{+\infty}$ and $d^{-\infty}$). Fix $k$, and consider 
the elements  $z_n={x_k}^{n+1}{y_k}^{-n}$. They belong to the coset $C$,
because $J/K$ is abelian, and their stable norm goes to infinity with $n$.
Therefore $C$ contains infinitely many
conjugacy classes.
 Let us now construct $c,d$ as above. Choose $a,b \in K$ generating a free group
of rank 2. We first explain how to get $c$. There is a problem only
if $ua^{\infty}=a^{-\infty}$ and $ub^{\infty}=b^{-\infty}$ . In that case
there exists integers $p,q$ with $ua^pu^{-1}=a^{-p}$ and  $ub^qu^{-1}=b^{-q}$.
 We take $c=a^pb^q$, noting that $ucu^{-1}=a^{-p}b^{-q}$ is different from
$c^{-1}=b^{-q}a^{-p}$. 

Once we have $c$, we choose $c^*\in K$ with $<c,c^*>$ free of rank 2, and we
obtain $d$ by applying the preceding argument using $c^*$ and $cc^*$ instead
of $a$ and $b$. The group $<c,d>$ is free of rank 2 because $d$ is a positive
word in $c^*$ and $cc^*$.
 
\subsection{Automorphisms of Gromov hyperbolic groups}
Let now $\phi$ be an automorphism of the Gromov hyperbolic group $G$ and 
let $\| \|$ denote the word metric with respect to some finite generating set for $G$. The automorphism $\phi$ is called hyperbolic if there is an integer
$m$ and a number $\lambda > 1$ such that, for all $g\in G$ we have 
$max(\|\phi^m(g)\|, \|\phi^{-m}(g)\|)\geq  \lambda\|g\|$.
For example a pseudo-Anosov homeomorphism of a closed surface
of genus larger then one induces a hyperbolic automorphism on the level of fundamental group. Also, an automorphism of finitely generated free group
with no nontrivial periodic conjugacy classes is hyperbolic.

\begin{lemma} \cite{bf} The mapping torus $M(\phi)$ of a hyperbolic automorphism
$\phi$ is Gromov hyperbolic group.  
\end{lemma}

\begin{theorem} The Reidemeister number $R(\phi)$ is infinite if group
$G$ is Gromov hyperbolic, non-elementary, and $\phi$ is hyperbolic automorphism.
\end{theorem}
{\sc Proof} The proof immediately folows from lemmas 2, 3 and 4    .

Let us consider an outer automorphism  $\Phi \in Out G$ corresponding
to automorphism of $\phi \in AutG$ and viewed as a collection of ordinary
automorphisms $\alpha \in Aut G$. We define $\alpha ,\beta \in \Phi$ to
be isogredient if $\beta=i_h\cdot\alpha\cdot {i_h}^{-1}$ for
some $h\in G$, with $i_h(g)=hgh^{-1}$.

\begin{lemma}\cite{LL} The set $S(\Phi)$ of isogredience classes is infinite if group
$G$ is Gromov hyperbolic, non-elementary, and $\Phi$ has finite order in the group $OutG$.
\end{lemma}

 {\sc Proof}
Let $J$ be the subgroup of $Aut G$ consisting of all automorphisms whose
image in $Out G$ is a power of $\Phi$. The exact sequence
$ {1}\rightarrow K \rightarrow J \rightarrow<\Phi> \rightarrow {1},$
 with $K=G/Center$ and $<\Phi>$ finite, shows that $J$ is hyperbolic,
non-elementary. The set of automorphisms $\alpha \in \Phi$ is a coset
of $J \mod K$. If $\alpha, \beta \in \Phi$ are isogredient they are
conjugate in $J$. The proof is therefore concluded by applying 
lemma 3  .

\begin{lemma}\cite{L}
 The set $S(\Phi)$ of isogredience classes is infinite if group
$G$ is a free group $F_n$ and $\Phi \in Out(F_n)$ fixes a nontrivial
conjugacy class.
\end{lemma}

\begin{theorem}\cite{LL}
For any $\Phi \in Out G$ , with $G$ Gromov hyperbolic, non-elementary,
the set $S(\Phi)$ of isogredience classes is infinite.
\end{theorem}
 {\sc Proof}
We describe main steps of the proof.
 By lemma 5, we may assume that $\Phi$ has infinite order.
By Paulin's theorem \cite{p} $\Phi$ preserves some $R$-tree $T$
with nontrivial minimal small action of $G$( recall that an action of $G$ is small if all ars stabilisers are virtually cyclic; the action of $G$ on $T$
is always irreducible( no global fixed point, no invariant line, no invariant end)). This means that there is an $R$-tree  $T$ equipped with
an isometric action of $G$ whose length function satisfies
$l\cdot \Phi=\lambda l$ for some $\lambda \geq 1$.\\
Step 1. Suppose $ \lambda =1$. Then $S(\Phi)$ is infinite.\\
Step 2. Suppose $\lambda >1$. Assume that arc stabilisers are finite, and
there exists $N_0\in N$ such that, for every $Q\in T$, the action
of $Stab Q$ on $\pi_o(T-{Q})$ has at most $N_0$ orbits.
Then  $S(\Phi)$ is infinite.\\
Step 3. If $\lambda >1$, then $T$ has finite arc stabilisers.
If $\lambda >1$ then from work of Bestwina-Feighn \cite{bf} it follows
that there exists $N_0\in N$ such that, for every $ Q\in T$, the action
of $Stab Q$ on $\pi_o(T-(Q))$ has at most $N_0$ orbits.

\begin{theorem} The Reidemeister number $R(\phi)$ is infinite if group
$G$ is Gromov hyperbolic, non-elementary, and $\phi$ is any automorphism
of $ G$.
\end{theorem}
{\sc Proof}
 By definition, the automorphisms
$\beta=i_m\cdot\alpha$ and $\gamma=i_n\cdot\alpha$ are isogredient
if and only if there exists $g\in G$ with
 $\gamma=i_g\cdot\beta\cdot {i_g}^{-1}$,or equivalently $n=gm\alpha(g^{-1})c$ with $c$ in center of $G$.
 So, the set $S(\Phi)$ of isogredience
classes of automorphisms representing $\Phi$ may be identified to the set
of twisted conjugacy classes of $G$ $\mod$ its center.
If $\phi$ is automorphism of finite order in $ AutG$, then the theorem 
immediately follows from lemma 5 .
If an automorphism  $\phi$ has infinite order in  $ AutG$ then theorem
follows from  theorem 2 .

\subsection{ The Co-Hopf property}

A group $G$ is called co-Hopf if every monomorphism of $G$
into itself is an isomorphism. It is fairly immediate to see
that a freely decomposable group is not co-Hopf.

\begin{lemma}(\cite{s})
Let $G$ be a non-elementary, torsion-free, Gromov hyperbolic
group. Then $G$ is co-Hopf if and only if  $G$ is freely
indecomposable.
\end{lemma}

\begin{theorem} The Reidemeister number $R(\phi)$ is infinite if group
$G$ is Gromov hyperbolic, non-elementary,torsion free, freely indecomposable
and $\phi$ is any monomorphism of $G$ into itself.
\end{theorem}
{\sc Proof}
 The proof follows from lemma 7 and theorem 3   .

\subsection{ Reduction to Injective Endomorphisms}
Let $G$ be a group and $\phi : G\rightarrow G$ an endomorphism.
We shall call an element $x\in G $ nilpotent if there is
an $n\in \N $ such that $\phi^n(x)=id$. Let $N$ be the set 
of all nilpotent elements of $G$.
 \begin{lemma}
The set $N$ is a normal subgroup of $G$. We have $ \phi(N) \subset N$ and 
$\phi^{-1}(N)=N$. Thus $\phi$ induces an endomorphism $[\phi/N](xN):=\phi(x)N$.
The endomorphism  $[\phi/N]: G/N \rightarrow G/N$ is injective,
and we have $R(\phi)=R([\phi/N])$.
\end{lemma}

{\sc Proof}

(i) Let $x\in N, g\in G$. Then for some $n\in \N$ we have $\phi^n(x)=id$.
Therefore $\phi^n(gxg^{-1})=\phi^n(gg^{-1})=id$. This shows that $gxg^{-1}\in N$
so $N$ is a normal subgroup of $G$.

(ii) Let  $x\in N$ and choose $n$ such that $\phi^n(x)=id$. Then  $\phi^{n-1}(\phi(x))=id$ so $\phi(x)\in N$. Therefore  $ \phi(N) \subset N$.

(iii) If  $\phi(x)\in N$ then there is an $n$ such that $\phi^n(\phi(x))=id$. Therefore $\phi^{n+1}(x)=id$ so  $x\in N$. This show that
 $\phi^{-1}(N)\subset N$. The converse inclusion follows from (ii).

(iv) We shal write $ \cR(\phi)$ for the set of $\phi$-conjugacy classes
of elements of $G$. We shall now  show that the map $x\rightarrow xN$
induces a bijection  $ \cR(\phi)\rightarrow \cR([\phi/N])$.
Suppose $x,y\in G$ are $\phi$-conjugate. Then there is a $g\in G$ with
$gx=y\phi(g).$ Projecting to the quotient group $G/N$ we have
$gnxN=yN\phi(g)N$, so $gNxN=yN[\phi/N](gN)$. This means that $xN$ and $yN$
are $[\phi/N]$-conjugate in $G/N$. Conversely suppose that $xN$ and $yN$
are  $[\phi/N]$-conjugate in $G/N$. Then there is a $gN \in G/N$
such that  $gNxN=yN[\phi/N](gN)$. In other words $gx\phi(g)^{-1}y^{-1})=id.$
Therefore $\phi^n(g)\phi^n(x)=\phi^n(y)\phi^n(\phi(g)).$

This shows that $\phi^n(x)$ and  $\phi^n(y)$ are $\phi$-conjugate.
Howeever by lemma 1 $x$ and  $\phi^n(x)$ are $\phi$-conjugate as are
 $y$ and  $\phi^n(y)$. Therefore $x$ and $y$ are  $\phi$-conjugate .

(v) We have shown that $x$ and $y$ are $\phi$-conjugate iff $xN$ and $yN$
are  $\phi/N$-conjugate. From this it follows that $x\rightarrow xN$
induces a  bijection from  $ \cR(\phi)$ to $ \cR([\phi/N])$.
Therefore $R(\phi)=R([\phi/N])$.

\begin{theorem} The Reidemeister number $R(\phi)$ is infinite if group
$G/N$ is Gromov hyperbolic, non-elementary,torsion free, freely indecomposable
and $\phi$ is any endomorphism of $G$ into itself.
\end{theorem}
{\sc Proof}
 The proof follows from lemma 8 and theorems 3  and 4    .

\section{Reidemeister number of continuous map}

Let now $X$ to be a connected, compact
polyhedron and $f:X\rightarrow X$ to be a continuous map.
Let a specific lifting $\tilde{f}:\tilde{X}\rightarrow\tilde{X}$ be chosen
as reference.
Then every lifting of $f$ can be written uniquely
as $\gamma\circ \tilde{f}$, with $\gamma\in\Gamma$.
So  the elements of $\Gamma$ serve as coordinates of
liftings with respect to the reference $\tilde{f}$.
Now,  for every $\gamma\in\Gamma$,  the composition $\tilde{f}\circ\gamma$
is a lifting of $f$;  so there is a unique $\gamma^\prime\in\Gamma$
such that $\gamma^\prime\circ\tilde{f}=\tilde{f}\circ\gamma$.
This correspondence $\gamma \rightarrow  \gamma^\prime$ is determined by
the reference $\tilde{f}$, and is obviously a homomorphism.
 
\begin{definition}
The endomorphism $\tilde{f}_*:\Gamma\rightarrow\Gamma$ determined
by the lifting $\tilde{f}$ of $f$ is defined by
$$
  \tilde{f}_*(\gamma)\circ\tilde{f} = \tilde{f}\circ\gamma.
$$
\end{definition}
 
It is well known that $\Gamma\cong\pi_1(X)$.
We shall identify $\pi=\pi_1(X,x_0)$ and $\Gamma$ in the following way.
Pick base points $x_0\in X$ and $\tilde{x}_0\in p^{-1}(x_0)\subset \tilde{X}$
once and for all.
 Now points of $\tilde{X}$ are in 1-1 correspondence with homotopy classes of paths
in $X$ which start at $x_0$:
for $\tilde{x}\in\tilde{X}$ take any path in $\tilde{X}$ from $\tilde{x}_0$ to 
$\tilde{x}$ and project it onto $X$;
conversely,  for a path $c$ starting at $x_0$, lift it to a path in $\tilde{X}$
which starts at $\tilde{x}_0$, and then take its endpoint.
In this way, we identify a point of $\tilde{X}$ with
a path class $\langle c \rangle $ in $X$ starting from $x_0$. Under this identification,
$\tilde{x}_0= \langle e \rangle $ is the unit element in $\pi_1(X,x_0)$.
The action of the loop class $\alpha =  \langle a \rangle \in\pi_1(X,x_0)$ on $\tilde{X}$
is then given by
$$
\alpha =  \langle a \rangle   :  \langle c \rangle \rightarrow \alpha \cdot  c = \langle a\cdot c\rangle.
$$
Now we have the following relationship between $\tilde{f}_*:\pi\rightarrow\pi$
and
$$
f_*  :  \pi_1(X,x_0) \longrightarrow \pi_1(X,f(x_0)).
$$

\begin{lemma}
Suppose $\tilde{f}(\tilde{x}_0) = \langle w \rangle$.
Then the following diagram commutes:
$$
\begin{array}{ccc}
  \pi_1(X,x_0)  &  \stackrel{f_*}{\longrightarrow}  &  \pi_1(X,f(x_0))  \\
                &  \tilde{f}_* \searrow \;\; &  \downarrow w_*   \\
                &                           &  \pi_1(X,x_0)
\end{array}
$$
where $w_*$ is the isomorphism induced by the path $w$.
\end{lemma}
In other words, for every  $\alpha = \langle a \rangle \in\pi_1(X,x_0)$, we have
$$
\tilde{f}_*(\langle a\rangle )=\langle w(f \circ a)w^{-1}\rangle.  
$$ 
\begin{remark}
In particular, if $ x_0 \in p(Fix( \tilde f))$ and $ \tilde x_0 \in Fix( \tilde f)$, then $ \tilde{f}_*=f_*$.
\end{remark} 
 
\begin{lemma}
Lifting classes of $f$ (and hence  fixed point classes, empty or not)
 are in 1-1 correspondence with
$\tilde{f}_*$-conjugacy classes in $\pi$,
the lifting class $[\gamma\circ\tilde{f}]$ corresponding
to the $\tilde{f}_*$-conjugacy class of $\gamma$.
We therefore have $R(f) = R(\tilde{f}_*)$.
\end{lemma}
 
We shall say that the fixed point class $p(Fix(\gamma\circ\tilde{f}))$,
which is labeled with the lifting class $[\gamma\circ\tilde{f}]$,
{\it corresponds} to the $\tilde{f}_*$-conjugacy class of $\gamma$.
Thus $\tilde{f}_*$-conjugacy classes in $\pi$ serve as coordinates for fixed
point classes of $f$, once a reference lifting $\tilde{f}$ is chosen.

\begin{theorem}
 Reidemeister number $R(f)$ is infinite for  any  homeomorphism $f$
of a compact polyhedron $X$ with the non-elementary, Gromov hyperbolic fundamental group $\pi_1(X,x_0)$.
\end{theorem}
 {\sc Proof}
Proof immediately follows from theorem 3 and lemma 10  .

\begin{theorem} 
The Reidemeister number $R(f)$ is infinite if fundamental group $\pi_1(X,x_0)$ is Gromov hyperbolic, non-elementary,torsion free, freely indecomposable and 
the map $f$ induces a monomorphism $\tilde{f}_*$  of $\pi_1(X,x_0)$  into itself.
\end{theorem}

{\sc Proof}
The proof follows from lemma 10  and  theorem 4   .

Let now $f$ is any continious map of $X$ .
Let $N$ be the set 
of all nilpotent elements of  $\pi_1(X,x_0)$ under  $\tilde{f}_*$ . 

\begin{theorem} 
 For any map  $f$ the Reidemeister number $R(f)$ is infinite  if the group $\pi_1(X,x_0)/N$ is Gromov hyperbolic, non-elementary,torsion free, freely indecomposable.
\end{theorem}

{\sc Proof}
The proof follows from lemma 10 and  theorem 5  .

\section{Pseudo-Anosov homeomorphism and asymptotic expansion}

Now we assume $X$ to be a compact surface of negative Euler characteristic and
 $f:X\rightarrow X$ is a pseudo-Anosov homeomorphism, i.e. there is a number $\lambda >1$ 
and a pair of transverse measured foliations $(F^s,\mu^s)$ and $(F^u,\mu^u)$ 
such that $f(F^s,\mu^s)=(F^s,\frac{1}{\lambda}\mu^s)$ and $f(F^u,\mu^u)=(F^u,\lambda\mu^u)$.
The mapping torus $T_f$ of $f:X\rightarrow X$ is the space obtained from $X\times [0 ,1]$ 
 by identifying $(x,1)$ with $(f(x),0)$ for all $x\in X$.
Take the base point $x_0$ of $X$ as the base point of $T_f$. According to van Kampen's Theorem, 
the fundamental group $\pi_ 1(T_f,x_0)$ is obtained from $\pi _1(X,x_0) $ by adding a new 
generator $z$ and adding the relations $z^{-1}gz=\tilde f_*(g)$ for all $g\in \pi _1(X,x_0)$, 
where $z$ is the generator of $\pi_1(S^1,x_0)$. 
This means that $\pi_ 1(T_f,x_0)$ is a semi-direct product of $\pi _1(X,x_0)$ with $Z$.
 
From Lemma 2 it follows that two elements $x,y$ of $\pi _1(X,x_0)$ are $\tilde f_*$-conjugate if and only if  $xz$ and $yz$ are conjugate 
in the usual sense in $\pi_ 1(T_f,x_0)$. Therefore $R(f) = R(\tilde{f}_*)$ is the number of usual conjugacy 
classes in the coset $\pi _1(X,x_0) \cdot z$ of $\pi _1(X,x_0)$ in $\pi_1(T_f,x_0)$.
There is a canonical projection
 $\tau : T_f \to R/Z$ given by
 $(x,s)\mapsto s$.
This induces a map
 $\pi_1(\tau):G=\pi_1(T_f,x_0)\to Z$.

We see that the Reidemeister number $R(f)$
 is equal to the number of homotopy
 classes of closed paths $\gamma$ in $T_f$
 whose projections onto $R/Z$
 are homotopic to the path
 $$
\begin{array}{cccl}
 \sigma :&[0,1] & \to     & $ $ \!\!\! R/Z \\
         &  s   & \mapsto & s.
\end{array}
 $$

Corresponding to this,  there is a group-theoretical interpretation of
 $R(f)$ as the number of
 usual conjugacy classes of elements $\gamma\in\pi_1(T_f)$
 satisfying $\pi_1(\tau)(\gamma) = z$.

 \begin{lemma}\cite{thu, otal}
The interior of  the mapping torus $ Int(T_f)$ admits a hyperbolic structure of finite volume if and only if
$f$ is isotopic to a pseudo-Anosov homeomorphism.
\end{lemma}

So, if the surface $X$ is closed and $f$ is isotopic to a pseudo-Anosov homeomorphism, 
the mapping torus $T_f$ can be realised as a hyperbolic 3-manifold ,  $H^3/\pi_ 1(T_f,x_0)$, where $H^3$
is the Poincare upper half space $\{(x,y,z): z  > 0, (x,y) \in R^2\}$ with the metric
$ds^2=(dx^2 +dy^2 + dz^2)/z^2$. 

The closed geodesics on a hyperbolic manifold are in one-to-one correspondence with
the free homotopy classes of loops. These classes of loops are in one-to-one correspondence with the 
conjugacy classes of loxodromic elements in the fundamental group of the hyperbolic manifold.

This correspondence allows  Ch. Epstein ( see \cite{eps}, p.127) to study the asymptotics of such functions 
as   $p_n(x)= \# \{$ primitive closed geodesics of length  less than $x$ represented by an element of the
form $gz^n \}$  using the Selberg trace formula. A primitive closed geodesic is one which is not 
an iterate of another closed geodesic. \\

  Phillips and Sarnak \cite{ps} have generalised results of Epstein in the following way. 
Let $M$ be a compact Riemannian manifold of negative curvature. Let $\psi : G  \rightarrow
 \Lambda $ be a surjective homomorphism of the fundamental group $G$ of $M$ onto an abelian
group $\Lambda$. Let $r$ be the rank of $\Lambda$ and let $m$ be the order of the torsion
of $\Lambda $; i.e.,  $\Lambda$ is isomorphic to $ Z^r \cdot F$,  $F$ being finite of order $m$.
For $\beta \in \Lambda$ let $p_{\beta}(x)$ be the number of primitive closed geodesics
$\gamma $ on $M$ of length at most $x$ for which $\psi (\gamma)=\beta$ .

\begin{theorem}\cite{ps}
We have the asymptotic expansion
$$p_{\beta}(x) \sim \frac{e^{(n-1)x}}{mx^{r/2 + 1}}( C_0 + \frac{C_1}{x} + \frac{C_2}{x^2}+ . . . .)$$
as $x \rightarrow \infty $. Here $C_0$ is independent of $\beta $; $C_0$ is the determinant of
a certain period matrix of harmonic 1-forms on $M$.
\end{theorem}

The proof of this result makes  routine use of the Selberg trace formula.
If $\Lambda = H_1(M, Z) $ then theorem 9 gives the asymptotic of the number of 
primitive closed geodesics of length at most $x$ lying in fixed homology class.
  In the more general case
of variable negative curvature similar  asymptotic was obtained by Pollicott and Sharp \cite{pos}.
They used dynamical approach based on the geodesic flow.
We will only need an asymptotic   for $p_1(x)$.  Note, that closed geodesics which are represented by an
element of the form $gz$ are  automatically primitive, because they wrap 
exactly once around the mapping torus (once around generator $z$). 
Taking $\Lambda=Z$ ,a surjective homomorphism  $\pi_1(\tau):G=\pi_1(T_f,x_0)\to Z$ and $\beta =z$-generator of 
group $Z$ in theorem 9 we have the  following asymptotic expansion \cite{eps, ps, pos, an} 
\begin{equation}
p_1(x)=  \frac{e^{hx}}{x^{3/2}}( \sum_{n=0}^{N}\frac{C_n}{x^{n/2}} + o(\frac{1}{x^{N/2}}) ),
\end{equation}
for any $N>0$,
where $h=\dim T_f  - 1=2$ is the topological entropy of the geodesic flow on the unit-tangent bundle 
$ST_f$, the constant  $C_0 > 0$ , constants $C_n$ vanish if $n$ is odd.
So in particular, using $N=1$ we have 
\begin{equation}
p_1(x) \sim  C_0 \frac{e^{hx}}{x^{3/2}},   x \rightarrow \infty
\end{equation}
Since $C_0 > 0$ we see that $p_1(x)$ goes to infinity as $x$ does.

 Now, using  one-to-one correspondences  in lemma 2  and   lemma 10  we  define  the  norm of  a fixed point class, or of the corresponding  
lifting class, or of the corresponding  twisted 
conjugacy class  $\{g\}_{\tilde f_*}$ in the fundamental group
of the surface $\pi _1(X, x_0)$,  as the length of the primitive closed geodesic  $\gamma $  on  $T_f$,
 which is represented by an element of the form $gz$. So, for example, the norm function $l^*$ 
on the set of twisted  conjugacy classes equals  $l^*=l\circ B$,
where $l$ is  length function on geodesics ($l(\gamma)$ is the length of the 
primitive closed geodesic  $\gamma$) and $B$ is a bijection
between the set of twisted 
conjugacy classes  $\{g\}_{\tilde f_*}$ in the fundamental group
of  the surface $\pi =\pi _1(X, x_0)$  and the set of closed geodesics represented by an
elements of the form $gz$ in the fundamental group $G :=\pi_ 1(T_f,x_0)$ .
We introduce the following  counting functions\\ 
 FPC($x$)= \# \{fixed point classes of $f$ of norm  less than $x$\},\\ 
L($x$)= \# \{lifting  classes of $f$ of norm  less than $x$\},\\
Tw($x$)= \#\{twisted conjugacy  classes for  $\tilde f_*$  in the fundamental group
of surface  of norm  less than $x$\}. 
We have following asymptotic expansion 
\begin{theorem}
Let $X$ be a closed surface of negative Euler characteristic and let $f:X\rightarrow X$ 
be a pseudo-Anosov homeomorphism. Then 
$$ FPC(x)=L(x)=Tw(x)= 
 \frac{e^{2x}}{x^{3/2}}( \sum_{n=0}^{N}\frac{C_n}{x^{n/2}}  + o(\frac{1}{x^{N/2}}) ),$$
where the constant  $C_0 > 0$ 
depends on the volume of  the  hyperbolic 3-manifold  $T_f$, and the
constants $C_n$ vanish if $n$ is odd.
\end{theorem}

{\sc Proof}
The proof follows from lemmas 2  and 10  and   asymptotic expansion (1).

\begin{corollary}
For pseudo-Anosov homeomorphisms  of  closed surfaces  the Reidemeister
number is infinite.
\end{corollary}
As we see, Theorem 10  gives   analytical proof of Lemma 3
in the case of pseudo-Anosov automorphism of surface group.

Theorem 10 indicates also, that we can develop an analog of Selberg theory
for twisted conjugacy classes of group automorphism.
We have the following question in this direction
\begin{question}
Is there an  asymptotic  of counting function 
for twisted conjugacy classes  of hyperbolic
automorphism of Gromov hyperbolic group?
\end{question}

\bigskip
 
 Institut f\"ur Mathematik, 

Ernst-Moritz-Arndt-Universit\"at Greifswald 
 
Jahn-strasse 15a, D-17489 Greifswald, Germany.
 
{\it E-mail address}: felshtyn@mail.uni-greifswald.de
\bigskip

\end{document}